\documentclass[11pt]{article}

\usepackage{amssymb}
\usepackage{amsmath}
\usepackage{theorem}
\usepackage{epsfig}
\usepackage{verbatim}
\usepackage{graphicx}

\usepackage{color}

\newtheorem{thm}{Theorem}[section]

\newtheorem{defi}[thm]{Definition}
\newtheorem{rem}[thm]{Remark}

\newcommand{\R}{\Bbb{R}}
\newcommand{\C}{\Bbb{C}}

\newcommand{\T}{\mathbb{T}}
\newcommand{\D}{\displaystyle}

\newcommand{\grad}{\nabla}

\newcommand{\gradp}{\grad^{\bot}}

\newcommand{\da}{\partial_\alpha}

\newcommand{\al}{\alpha}

\newcommand{\ztil}{\tilde{z}}

\newcommand{\pa}{\partial}

\newenvironment{sket}{\begin{trivlist} \item[] {\em Sketch of the Proof:}}{\hfill $\Box$
                       \end{trivlist}}

\begin{document}

\title{Splash singularity for water waves}
\date{}
\author{ Angel Castro, Diego C\'ordoba, Charles Fefferman, \\ Francisco Gancedo and Javier G\'omez-Serrano}


\title{Splash singularity for water waves}
\date{ }


\maketitle

\begin{abstract}

We exhibit smooth initial data for the 2D water wave equation for which we prove that smoothness of the interface breaks down in finite time. Moreover, we show a stability result together with numerical evidence that there exist solutions of the 2D water wave equation that start from a graph, turn over and collapse in a splash singularity (self intersecting curve in one point) in finite time.


\end{abstract}


\section{Introduction}

\begin{figure}[h!]\centering
\includegraphics[scale=0.2]{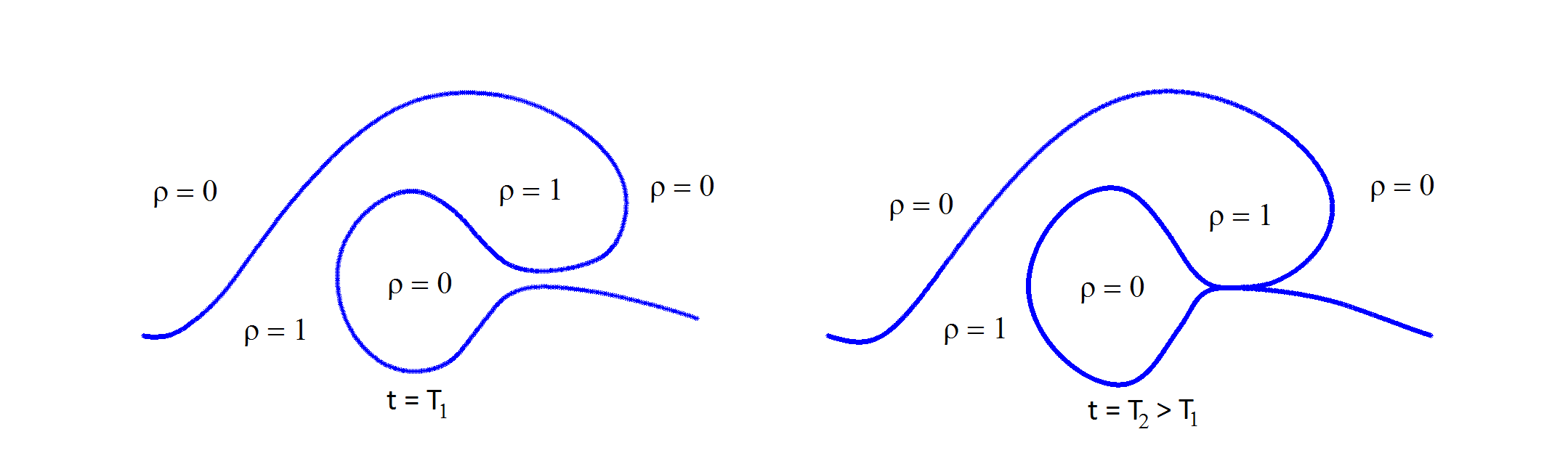}
\caption{Splash singularity. The interface collapses at one point.}
\label{qualitativesplash}
\end{figure}

We consider the 2D water wave equation, which governs the motion of the interface between a 2D inviscid incompressible irrotational fluid and a vacuum, taking gravity into account but neglecting surface tension. We prove that an initially smooth interface may in finite time become singular by the mechanism illustrated in fig. \ref{qualitativesplash}. We call such  a singularity a ``splash''. We also present numerical evidence for a scenario in which the interface starts out as  a smooth graph, then ``turns over" after finite time, and finally produces a splash, as in fig. \ref{PictureNonTildaZoom}.

The equations of motion in $\mathbb{R}^{2}$ for the density $\rho = \rho(x,t)$, $x\in\mathbb{R}^2$, $t \geq 0$, the velocity $v = (v^{1},v^{2})$ and the pressure $p = p(x,t)$ are:

\begin{equation}
\label{euler}
\left\{
\begin{array}{ccc}
\rho(v_t+v\cdot\grad v) & = &-\nabla p-(0,\rho), \\
\rho_t + v\cdot\nabla\rho &= &0, \\
\nabla\cdot v &= &0.
\end{array}
\right.
\end{equation}
Above, the acceleration due to gravity is taken equal to one for the sake of simplicity.

The free boundary is parameterized by
$$
\partial \Omega^j(t)=\{z(\al,t)=(z_1(\al,t),z_2(\al,t)):\al\in\R\},
$$
where the regions $\Omega^{j}(t)$ are defined by

\begin{equation*}\label{density}
\rho(x_1,x_2,t)=\left\{\begin{array}{cl}
                    0,& x\in\Omega^1(t)\\
                    1,& x\in\Omega^2(t)=\mathbb{R}^2 - \Omega^1(t).
                 \end{array}\right.
\end{equation*}

We assume that the fluid is irrotational, i.e. the vorticity $\gradp\cdot v = 0$, in the interior of each domain $\Omega^j$
($j=1,2$). The vorticity will be supported on the free boundary curve
$z(\al,t)$ and it has the form
\begin{equation*}
\gradp\cdot v(x,t)=\omega(\al,t)\delta(x-z(\al,t)),
\end{equation*}
i.e. the vorticity is a Dirac measure defined by
$$
\langle\gradp\cdot v,\eta\rangle=\int_{\R}\omega(\al,t)\eta(z(\al,t))d\al,
$$
with $\eta(x)$ a test function. We present results for the following geometries:
\begin{enumerate}
\item open curves asymptotic to the horizontal at infinity
$$\D\lim_{\al\rightarrow\infty}(z(\al,t)-(\al,0))=0.$$
\item periodic curves in the space variable
\begin{equation*}
z(\al+2k\pi,t)=z(\al,t)+2k\pi(1,0).
\end{equation*}
\item closed contours $$z(\al+2k\pi,t) = z(\al,t).$$
\end{enumerate}

However, in this paper we will only deal with the first case as the others are similar. One method to derive the equations for the evolution of $v$, $\rho$ is to write the velocity as the orthogonal gradient of the stream function, take the curl and recover the velocity by inverting the Laplacian, i.e.,  we apply the Biot-Savart law. Here we use the fact that the vorticity is concentrated on the interface;
\begin{equation*}\label{BS}
v(x,t)=\gradp\Delta^{-1}(\gradp\cdot v)(x,t)=\frac{1}{2\pi}\int_\R
\frac{(x-z(\al,t))^{\bot}}{|x-z(\al,t)|^2}\omega(\al,t)d\al,
\end{equation*}
with $x\neq z(\al,t)$.

Taking limits of the above equation, by approaching the boundary in the normal direction, we obtain the velocity of the interface, to which we can add any term $c$ in the tangential direction $z_{\al}$ without
modifying the geometry of the interface. Thus the interface satisfies
\begin{equation}\label{em}
z_t(\al,t)=BR(z,\omega)(\al,t)+c(\al,t) z_{\al}(\al,t),
\end{equation}
where the Birkhoff-Rott integral is defined by
\begin{equation*}\label{BR}
BR(z,\omega)=\frac{1}{2\pi}P.V.\int_\R
\frac{(z(\al,t) - z(\beta,t))^{\bot}}{|z(\al,t) - z(\beta,t)|^2}\omega(\al,t)d\al.
\end{equation*}

The system is closed by using Euler equations (for details see for example \cite{Cordoba-Cordoba-Gancedo:interface-water-waves-2d});

\begin{align}
\begin{split}\label{cEuler}
\omega_t(\al,t)&=-2\partial_t BR(z,\omega)(\al,t)\cdot
\da z(\al,t)-\da( \frac{|\omega|^2}{4|\da z|^2})(\al,t) +\da (c \omega)(\al,t)\\
&\quad +2c(\al,t)\da BR(z,\omega)(\al,t)\cdot\da z(\al,t)+2\da
z_2(\al,t).
\end{split}
\end{align}
Then the dynamic equations,  for the interface  $z=z(\al,t)$ and the  vorticity are the system given by \eqref{em} and \eqref{cEuler} and are known as the water-wave equations. 

Taking the divergence of the Euler equation \eqref{euler} and recalling that the flow is irrotational in the interior of the regions $\Omega^{j}(t)$, we find that
$$ -\Delta p = |\nabla v|^{2} \geq 0$$
which, together with the fact that the pressure is zero on the interface implies by Hopf's lemma in $\Omega^{2}(t)$ that $$\sigma(\al,t) \equiv -|z_{\al}^{\perp}(\al,t)|\partial_{n} p(z(\al,t),t) > 0,$$ where $\partial_{n}$ denotes the normal derivative. This is known as the Rayleigh-Taylor condition (see \cite{Rayleigh:instability-jets} and \cite{Taylor:instability-liquid-surfaces-I}) which was first proved by Wu in \cite{Wu:well-posedness-water-waves-2d} and \cite{Wu:well-posedness-water-waves-3d}.

 The first results concerning the Cauchy problem for the linearized
version of water waves and small data in Sobolev spaces are due to Craig \cite{Craig:existence-theory-water-waves}, Nalimov \cite{Nalimov:cauchy-poisson}, Beale et al. \cite{Beale-Hou-Lowengrub:growth-rates-linearized}
and Yoshihara \cite{Yosihara:gravity-waves}. The well-posedness in Sobolev spaces for the water-wave problem was proven by Wu in \cite{Wu:well-posedness-water-waves-2d} with the assumption that the initial free interface is non-self-intersecting (satisfies the arc-chord condition). For  recent work on local existence  see Wu \cite{Wu:well-posedness-water-waves-3d}, Christodoulou-Lindblad \cite{Christodoulou-Lindblad:motion-free-surface}, Lindblad
  \cite{Lindblad:well-posedness-motion},  Coutand-Shkoller \cite{Coutand-Shkoller:well-posedness-free-surface-incompressible}, Shatah-Zeng  \cite{Shatah-Zeng:geometry-priori-estimates}, Zhang-Zhang
    \cite{Zhang-Zhang:free-boundary-3d-euler}, C\'ordoba et al. \cite{Cordoba-Cordoba-Gancedo:interface-water-waves-2d}, Lannes \cite{Lannes:well-posedness-water-waves}, \cite{Lannes:stability-criterion}, Alazard-Metivier \cite{Alazard-Metivier:paralinearization} and Ambrose-Masmoudi \cite{Ambrose-Masmoudi:zero-surface-tension-waterwaves}.  The issue of long time existence has been treated in Alvarez-Lannes \cite{AlvarezSamaniego-Lannes:large-time-existence-water-waves} where well-posedness over large time scales is shown and different asymptotic regimes are justified. Wu proved in \cite{Wu:almost-global-wellposedness-2d} exponential time of existence in two dimensions for small initial data, and Germain et al. in \cite{Germain-Masmoudi-Shatah:global-solutions-gravity-water-waves} and Wu in \cite{Wu:global-wellposedness-3d} global existence for small data in the three dimensional case  (two dimensional interface). In \cite{Castro-Cordoba-Fefferman-Gancedo-LopezFernandez:rayleigh-taylor-breakdown} and \cite{Castro-Cordoba-Fefferman-Gancedo-LopezFernandez:turning-waves}, Castro et al. showed that there exist large initial data parameterized as a graph for which in finite time the interface reaches a regime in which it is no longer a graph. For previous numerical simulations showing this phenomenon see Beale et al. \cite{Beale-Hou-Lowengrub:convergence-boundary-integral}.

The outline of the paper is the following: in section 2 we will describe the equations in a transformed domain which will circumvent the problem of having a singularity where the arc-chord condition fails as the curve self-intersects, i.e. a splash singularity forms. In section 3 we will outline the proof of a local existence theorem for the equations in the new domain, both for the analytic case and for Sobolev spaces. Section 4 will be devoted to a stability theorem, whereas section 5 will comment on the numerical results obtained towards the splash singularity starting from a graph. Finally in section 6 we describe the ideas that we hope will lead to  a computer assisted proof of the existence of a solution that starts as a graph and ends in a splash.

\section{The equations in the tilde domain}

In this section we will rewrite the equations by applying a transformation from the original coordinates to new ones which we will denote with a tilde. The purpose of this transformation is to be able to deal with the failure of the arc-chord condition. We start by reformulating the set of equations, in the non tilde domain, for the case of a periodic contour in terms of the velocity potential. From \eqref{euler} and since $v$ is irrotational in $\Omega^2(t)$ we have that:

\begin{align}
\label{nontildapotential}
\nabla \times v & = 0 \quad \text{in $\Omega^2(t)$}\nonumber\\
\nabla \cdot v & = 0\quad \text{in $\Omega^2(t)$} \nonumber\\
\phi_t(x,y,t) + \frac{1}{2}|v(x,y,t)|^{2} & = -(p(x,y,t) - p^{*}(t)) - y\quad \text{in $\Omega^2(t)$}\nonumber\\
z_t(\al,t)&=u(\alpha,t)+c(\al,t)z_\alpha(\alpha,t)\nonumber\\
\Phi_\alpha(\alpha,t) & = u(\alpha,t)\cdot z_\alpha(\alpha,t)\nonumber\\
p&=\text{constant}\equiv0\quad \text{in $\Omega^{1}(t)$}\nonumber\\
z(\alpha,0)&=z^0(\alpha)\nonumber\\
\Phi_\alpha(\alpha,0)&=\Phi_\alpha^0(\alpha)
\end{align}
where $\phi$ is the velocity potential, $\Phi(\alpha,t)$ is its limit at the interface coming from the fluid region, $\nabla \phi = v$, $p^{*}$ is a function of $t$ alone, $c$ is a free quantity which represents the reparameterization freedom,  $u(\al,t)$ is the limit of the velocity at the interface coming from the fluid region and $\Phi_\alpha^0(\alpha)$ has zero mean. We also want  the velocity to be in $L^2\left(\Omega^2(t) \cap \left([-\pi,\pi]\times \R\right)\right)$ and periodic in the $x-$coordinate i.e.
$$v(x+2\pi,y)=v(x,y)\quad \text{in $\Omega^2(t)$}.$$

Note that $\Phi$ is periodic in the horizontal variable, because $v$ is periodic and $v(x,y)$ tends to zero as $y$ tends to $-\infty$.

In order to simplify this system we  use the stream function $\psi$ and consider the equations

\begin{align}
\label{nontildastream}
\Delta \psi(x,y,t) & = 0 \quad \text{in $\Omega^2(t)$} \nonumber\\
\left.\partial_{n} \psi\right|_{z(\al,t)} & =- \frac{\Phi_{\al}(\al,t)}{|z_{\al}(\al,t)|} \nonumber\\
\psi(x+2\pi,y,t)&=\psi(x,y,t)\quad \text{in $\Omega^2(t)$}\nonumber\\
\psi(x,y) & \text{ is $O(1)$ as $y \to -\infty$}\nonumber\\
v&\equiv\nabla^\perp \psi\quad \text{in $\Omega^2(t)$}\nonumber\\
z_t(\al,t) & = u(\al,t) + c(\al,t)z_{\al}(\al,t) \nonumber\\
\phi(x,y,t) & \text{ is the harmonic conjugate of $\psi(x,y,t)$ in $\Omega^2$}\nonumber\\
\Phi_t(\al,t) & = \frac{1}{2}|u(\al,t)|^2 + c(\al,t)u(\al,t) \cdot z_{\al}(\al,t) - z_{2}(\al,t) + p^*(t) \nonumber\\
z(\alpha,0)&=z^0(\alpha)\nonumber\\
\Phi_\alpha(\alpha,0)&=\Phi_\alpha^0(\alpha)
\end{align}

Although we may take as an initial condition the tangential component of the velocity multiplied by the modulus of the tangent vector, i.e. $\Phi_\alpha$, we can also solve the system \eqref{nontildastream} by taking as an initial condition the normal component of the initial velocity multiplied by the modulus of the normal vector, i.e. $\Psi_\alpha$, ($\Psi(\alpha,t)=\psi(z(\alpha,t),t)$), as we can transform one into the other.

It can be checked that solutions of the system \eqref{nontildastream} are also solutions of the system \eqref{nontildapotential}.
Let us consider $\ztil(\al,t) = (\tilde{z}_1(\al,t), \tilde{z}_2(\al,t)) \equiv P(z(\al,t))$ where $P$ is a conformal map defined in the water region that will be given as:
$$ P(w) = \left(\tan\left(\frac{w}{2}\right)\right)^{1/2},\quad w\in\C,$$
 for a branch of the square root that separates the self-intersecting points of the interface. Here $P(w)$ will refer to a 2 dimensional vector whose components are the real and imaginary parts of $P(w_1 + iw_2)$.  In this setting, $P^{-1}(z)$ will be well defined modulo multiples of $2\pi$.

The water wave equations are invariant under time reversal. To obtain a solution that ends in a splash, we can therefore take our initial condition to be a splash, and show that there is a smooth solution for small times $t > 0$. As initial data we are interested in considering a curve that intersects itself at one point, as in fig. \ref{PictureNonTildaZoom}. More precisely, we will use as initial data \emph{splash curves} which are defined as follows:
\begin{defi}
We say that $z(\al)$ is a \emph{splash curve} if
\begin{enumerate}
\item $z_{1}(\al) - \al, z_2(\al)$ are smooth functions and $2\pi$-periodic.
\item $z(\al)$ satisfies the arc-chord condition at every point except at $\alpha_1$ and $\alpha_2$, with $\alpha_1 < \alpha_2$ where $z(\al_1) = z(\al_2)$ and $|z_{\al}(\al_1)|, |z_{\al}(\al_2)| > 0$. This means $z(\al_1) = z(\al_2)$, but if we remove either a neighborhood of $\al_1$ or a neighborhood of $\al_2$ in parameter space, then the arc-chord condition holds.
\item The curve $z(\alpha)$ separates the complex plane into two region; a  connected water region and a vacuum region. The water region contains each point $x+iy$ for which y is large negative. We choose the parametrization such that the normal vector $n=\frac{(-\pa_\alpha z^2(\alpha), \pa_\alpha z^1(\alpha))}{|\pa_\alpha z(\alpha)|}$ points to the vacuum region.
\item We can choose a branch of the function $P$ on the water region such that the curve $\tilde{z}(\al) = P(z(\al))$ satisfies:
\begin{enumerate}
\item $\tilde{z}_1(\al)$ and $\tilde{z}_2(\al)$ are smooth and $2\pi$-periodic.
\item $\tilde{z}$ is a closed contour.
\item $\tilde{z}$ satisfies the arc-chord condition.
\end{enumerate}
We will choose the branch of the root that produces that
$$ \lim_{y \to -\infty}P(x+iy) = -e^{-i \pi/4}$$
independently of $x$.
\item $P(w)$ is analytic in $w$ and $\frac{dP}{dw}(w) \neq 0$ if $w$ belongs to the water region.
\item $\tilde{z}(\al) \neq q^l$ for $l=0,...,4$, where
\begin{equation}\label{points}
q^0=\left(0,0\right),\quad
q^1=\left(\frac{1}{\sqrt{2}},\frac{1}{\sqrt{2}}\right),\quad
q^2=\left(\frac{-1}{\sqrt{2}},\frac{1}{\sqrt{2}}\right),\quad
q^3=\left(\frac{-1}{\sqrt{2}}, \frac{-1}{\sqrt{2}}\right),\quad
q^4=\left(\frac{1}{\sqrt{2}}, \frac{-1}{\sqrt{2}}\right).
\end{equation}
\end{enumerate}
\end{defi}

\begin{figure}[h!]\centering
\includegraphics[scale=0.3]{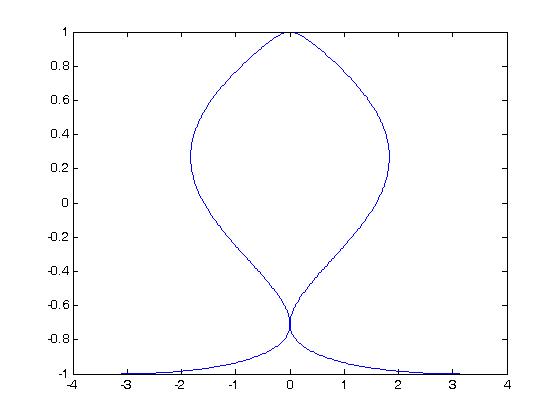}
\includegraphics[scale=0.3]{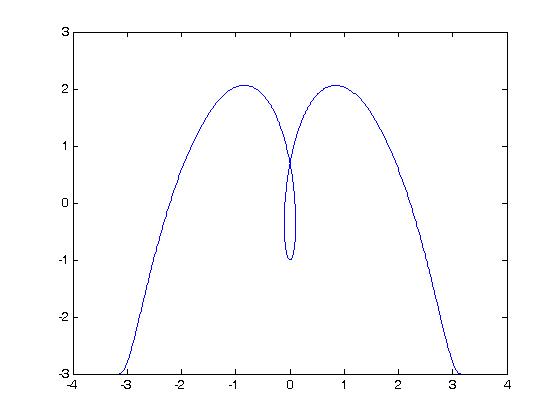}
\caption{Two examples of  non-splash curves.}
\label{nosplash}
\end{figure}

From now on, we will always work with splash curves as initial data. Condition 6 will be used in the local existence theorems and can be proved to hold for short enough time as long as the initial condition satisfies it. We will also need that the interface passes below the points $(\pm\pi,0)$ (or, equivalently, that those points belong to the vacuum region) in order for the tilde region to be a closed curve and the vacuum region to lie on the outer part. For a splash curve this is trivial from the definition.  For more information about the transformation of both regions, check the figures \ref{PictureSplash} and \ref{PictureNonTildaZoom} and notice that we rule out the scenarios in fig. \ref{nosplash}.

We will now write \eqref{nontildastream} in the new tilde coordinates. We define the following quantities:

$$ \tilde{\psi}(x,y,t) \equiv \psi(P^{-1}(x,y,t)), \quad \tilde{\phi}(x,y,t) \equiv \phi(P^{-1}(x,y,t)), \quad \tilde{v}(x,y,t) \equiv \nabla \tilde{\phi}(x,y,t).$$

Let us note that as $\psi$ and $\phi$ are $2\pi$ periodic, the resulting $\tilde{\psi}$ and $\tilde{\phi}$ are well defined. We do not have problems with the harmonicity of $\tilde{\psi}$ or $\tilde{\phi}$ at the point which is mapped from minus infinity (which belongs to the water region) by $P$ as $\phi$ and $\psi$ are well defined at infinity. Also, the periodicity of $\phi$ and $\psi$ causes $\tilde{\phi}$ and $\tilde{\psi}$ to be continuous (and harmonic) at the interior of $P(\Omega^2(t))$.

Let us assume that there exists a solution of \eqref{nontildastream} and that we take $u_n = \frac{\Psi_{\al}}{|z_{\al}|}$ such that $u_n(\al_1), u_n(\al_2) < 0$ for all $0 < t < T, T$ small enough, thus $z(\al,t)$ satisfies the arc-chord condition and does not touch the removed branch from $P(w)$.

Thus the system \eqref{nontildastream} in the new coordinates reads
\begin{align}
\label{tildastream}
\Delta \tilde{\psi}(x,y,t) & = 0 \quad \text{in $P(\Omega^2(t))$} \nonumber\\
\left.\partial_{n} \tilde{\psi}\right|_{\tilde{z}(\al,t)} & =- \frac{\tilde{\Phi}_{\al}(\al,t)}{|\tilde{z}_{\al}(\al,t)|} \nonumber\\
\tilde{v}&\equiv\nabla^\perp \tilde{\psi}\quad \text{in $P(\Omega^2(t))$}\nonumber\\
\tilde{z}_{t}(\al,t) & = Q^2(\al,t)\tilde{u}(\al,t) + c(\al,t)\tilde{z}_{\al}(\al,t) \nonumber\\
\tilde{\Phi}_{t}(\al,t) & = \frac{1}{2}Q^2(\al,t)|\tilde{u}(\al,t)|^2 + c(\al,t)\tilde{u}(\al,t) \cdot \tilde{z}_{\al}(\al,t) - P^{-1}_2(\tilde{z}(\al,t)) + p^{*}(t)\nonumber\\
\tilde{z}(\alpha,0)&=\tilde{z}^0(\alpha)\nonumber\\
\tilde{\Phi}_\alpha(\alpha,0)&=\tilde{\Phi}_\alpha^0(\alpha) = \Phi_{\al}^{0}(\al).
\end{align}
where $\tilde{u}$ is the limit of the velocity coming from the fluid region in the tilde domain and $Q^2(\al,t) = \left|\frac{dP}{dw}(z(\al,t))\right|^{2}$.

We can therefore solve the Neumann problem for the stream function
in the fluid domain with boundary conditions on the interface. In fact there exists a function $\tilde{\omega}$ satisfying
$$ \psi(\tilde{x},t) = \frac{1}{2\pi}\int_{-\pi}^{\pi}\log(|\tilde{x} - \tilde{z}(\al,t)|)\tilde{\omega}(\alpha,t)d\alpha, \quad \tilde{x} \text{ belonging to the fluid region, }$$
which implies
$$ \tilde{v}(\tilde{x},t) = \nabla^{\perp} \tilde{\psi}(\tilde{x},t) = \frac{1}{2\pi}P.V\int_{-\pi}^{\pi}\frac{(\tilde{x} - \tilde{z}(\al,t))^{\perp}}{|\tilde{x} - \tilde{z}(\al,t)|^2}\tilde{\omega}(\al,t)d\al.$$
Taking limits from the fluid region we obtain
$$ \tilde{u}(\al,t) = BR(\tilde{z},\tilde{\omega}) + \frac{\tilde{\omega}}{2|\tilde{z}_{\al}|^2}\tilde{z}_{\al}.$$
The evolution of $\tilde{\omega}$ is calculated in the following way. First, let us recall the equations
\begin{align}
\label{equations-hou-tilde}
\tilde{z}_{t}(\al,t) & = Q^2(\al,t)\tilde{u}(\al,t) + c(\al,t)\tilde{z}_{\al}(\al,t) \nonumber\\
\tilde{\Phi}_{t}(\al,t) & = \frac{1}{2}Q^2(\al,t)|\tilde{u}(\al,t)|^2 + c(\al,t)\tilde{u}(\al,t) \cdot \tilde{z}_{\al}(\al,t) - P^{-1}_2(\tilde{z}(\al,t)) \nonumber\\
\tilde{\Phi}_{\al}(\al,t) & = \tilde{u}(\al,t) \cdot \tilde{z}_{\al}(\al,t) \nonumber \\
\tilde{z}(\alpha,0)&=\tilde{z}^0(\alpha)\nonumber\\
\tilde{\Phi}_\alpha(\alpha,0)&=\tilde{\Phi}_\alpha^0(\alpha) = \Phi_{\al}^{0}(\al)
\end{align}

Those equations will be the ones used in section 5 with $c(\al,t) = 0$. Substituting the expression for $\tilde{u}(\al,t)$ and performing the change $\tilde{c}(\al,t) = c(\al,t) + \frac{1}{2}Q^2(\al,t)\frac{\tilde{\omega}(\al,t)}{|\tilde{z}_{\al}(\al,t)|^{2}}$ we obtain

\begin{align}\label{zeq}
\tilde{z}_{t}(\al,t) & = Q^2(\al,t)BR(\tilde{z},\tilde{\omega})(\al,t) + \tilde{c}(\al,t)\tilde{z}_{\al}(\al,t)
\end{align}

and the evolution equation for $\tilde{\omega}$
\begin{align}\label{eqomega}
\tilde{\omega}_{t}(\al,t) & = -2 \partial _{t} BR(\tilde{z},\tilde{\omega})(\al,t) \cdot \tilde{z}_{\al}(\al,t) - |BR(\tilde{z},\tilde{\omega})|^{2} \partial_{\al}Q^{2}(\al,t) \nonumber \\
& - \partial_{\al}\left(\frac{Q^2(\al,t)}{4}\frac{\tilde{\omega}(\al,t)^2}{|\tilde{z}_{\al}(\al,t)|^{2}}\right)
 + 2\tilde{c}(\al,t) \partial_{\al}BR(\tilde{z},\tilde{\omega}) \cdot \tilde{z}_{\al}(\al,t) \nonumber \\
& + \partial_{\al}\left(\tilde{c}(\al,t)\tilde{\omega}(\al,t)\right)
- 2\partial_{\al}\left(P^{-1}_2(\tilde{z}(\al,t))\right).
\end{align}

\begin{rem}
Equation \eqref{eqomega} is analogous to \eqref{cEuler}. In fact, if we  set $Q \equiv 1$ in \eqref{eqomega} we recover \eqref{cEuler}.
\end{rem}

Note that for the tilde domain, the Rayleigh-Taylor condition is the same as in the first domain, i.e:

$$ \nabla p(\al,t) \cdot z_{\al}^{\perp}(\al,t) = \nabla \tilde{p}(\al,t) \cdot \tilde{z}_{\al}^{\perp}(\al,t)$$

where $\tilde{p} = p \circ P^{-1}$.

Our strategy will be the following: we will consider the evolution of the solutions in the tilde domain and then see that everything works fine in the original domain.

We will have to obtain the normal velocity once given the tangential velocity, and viceversa. To do this, we just have to notice that
$$ \tilde{\Phi}_{\al}(\al,t)  = \tilde{u}(\al,t) \cdot \tilde{z}_{\al}(\al,t) =
BR(\tilde{z},\tilde{\omega})\cdot \tilde{z}_{\al}(\al,t) + \frac{\tilde{\omega}(\al,t)}{2}.$$
From that, we can invert the equation and get $\tilde{\omega}$, plug it into the following expression for $\tilde{\psi}$:
 $$ \tilde{\psi}(\tilde{x}) = \frac{1}{2\pi}\int_{-\pi}^{\pi}\log(|\tilde{x} - \tilde{z}(\al)|)\tilde{\omega}(\alpha)d\alpha$$
and restrict ourselves to the interface to get $\tilde{\Psi}(\al,t)$. Taking a derivative in $\al$ we can recover the normal component of the velocity. An analogous reasoning can be done to get the tangential velocity from the normal by solving the complementary Neumann problem for $\phi$.

We  now note that a solution of the system \eqref{tildastream} in the tilde domain gives rise to a solution of the system \eqref{nontildastream} in the non-tilde domain, by inverting the map $P$. In fact, this will be the implication used in Theorem \ref{localexistencenontilde} (finding a solution in the tilde domain, and therefore in the non-tilde).

\begin{rem}
A similar argument works for the other two settings (closed contour and asymptotic to horizontal) by choosing an appropriate $P(w)$ that separates the singularity.
\end{rem}


\section{Local existence at the splash}

The main result in this section is a local existence proof for the splash singularity. To avoid the arc-chord condition failure, we will prove the local existence in the tilde domain. This can be done in two different settings, namely in the space of analytic functions and the Sobolev space $H^{s}$.
\begin{thm}
\label{localexistencenontilde}
Let $z^{0}(\al)$ be a splash curve such that $z_{1}^{0}(\al) - \al, z_2^{0}(\al) \in H^{4}(\T)$.  Let $u^{0}(\al)\cdot (z^0_\alpha)^\perp(\alpha)\in H^4(\T)$ satisfying:
\begin{enumerate}
\item $\displaystyle \left(u^{0}\cdot\frac{(z_\al^0)^\bot}{|z_\al^0|}\right)(\al_1) < 0, \left(u^{0}\cdot\frac{(z_\al^0)^\bot}{|z_\al^0|}\right)(\alpha_2)<0,$
\item $\displaystyle \int_{\partial \Omega}u^{0}\cdot\frac{(z_\al^0)^\bot}{|z_\al^0|} ds = \int_{\T} u^{0}(\al)\cdot (z_\al^0)^\bot d\al = 0$.
\end{enumerate}
Then there exist a finite time $T > 0$, a time-varying curve $z(\al,t) \in C([0,T];H^{4}(\T))$ satisfying:
\begin{enumerate}
\item $z_{1}(\al,t) - \al, z_{2}(\al,t)$ are $2\pi$-periodic,
\item $z(\al,t)$ satisfies the arc-chord condition for all $t \in (0,T]$,
\end{enumerate}
and $u(\al,t)\in C([0,T];H^3(\T))$ which provides a solution of the water wave equations \eqref{nontildastream} with $z(\al,0)=z^0(\al)$ and $u(\al,0)\cdot(z_\alpha)^\perp(\alpha,0)=u^0(\al)\cdot(z^0_\alpha)^\perp(\alpha)$.
\end{thm}
\emph{Sketch of the proof:}
Using the fact that there is local existence to the initial data in the tilde domain and applying $P^{-1}$ to the solution obtained there, we can get a curve $z(\al,t)$ that solves the water waves equation in the non tilde domain. This leads to the proof of Theorem \ref{localexistencenontilde}. Details on the local existence in the tilde domain are shown below.


\subsection{Local existence for analytic initial data in the tilde domain}

In this subsection, we will work on the tilde domain, and all tildes will be dropped for the sake of simplicity.

We will work with $c = 0$. We have the following system:
\begin{equation}
\label{1paco}
\left\{
\begin{array}{rcl}
z_t & = & \left|\frac{dP}{dw}(P^{-1}(z))\right|^{2}u \\
\Phi_t & = & \frac{1}{2}\left|\frac{dP}{dw}(P^{-1}(z))\right|^{2}|u|^{2} - P_2^{-1}(z)\\
u & = & BR(z,\omega) + \frac{\omega}{2|z_{\al}|^{2}}z_{\al} \\
\Phi_{\al} & = & \frac{\omega}{2} + BR(z,\omega) \cdot z_{\al} \\
\left|\frac{dP}{dw}(P^{-1}(z(\al,t)))\right|^{2} & = & \frac{1}{16}\left|\frac{1+ (z_1(\al,t)+iz_2(\al,t))^{4}}{z_1(\al,t)+iz_2(\al,t)} \right|^2 \\
P_2^{-1}(z(\al,t)) & = & \ln\left|\frac{i+(z_1(\al,t)+iz_2(\al,t))^2}{i-(z_1(\al,t)+iz_2(\al,t))^2}\right|
\end{array}
\right.
\end{equation}


We demand that $z^0(\al) \neq (0,0)$ to find the function $\frac{dP}{dw}(P^{-1}(z(\al,t)))$ well defined. This condition is going to remain true for short time. We also consider $z^0(\al) \neq q^l$, $l=1,...,4$ in \eqref{points} to get $P_2^{-1}(z(\al,t))$ well defined. Again this is going to remain true for short time.

We consider the space
\begin{align*}
H^3(\partial S_r) = \left\{f \right.&\text{ analytic in } S_{r}= \{\al + i \eta, |\eta| < r\}:\, \|f\|_{r}^{2}<\infty  \\
&  \|f\|_r^2=\|f\|^2_{L^2(\partial S_r)}+\|\da^3f\|^2_{L^2(\partial S_r)} \text{ where }\\
& \left. \|f\|^2_{L^2(\partial S_r)}=\sum_{\pm}\int_{-\pi}^{\pi}|f(\al\pm ir)|^{2}d\al, f \; 2\pi \text{-periodic}\right\}
\end{align*}
and $(z,\Phi) \in (H^3(\partial S_r))^3 \equiv X_{r}$.

We have the following theorem:
\begin{thm}
\label{localexistencetilde}
Let $z^{0}(\al)$ be a splash curve and let $u^{0}\cdot\frac{z^0_\alpha}{|z^0_\alpha|}(\alpha)=\frac{\Phi^0_\alpha}{|z^0_\alpha|}(\alpha)$ be the initial tangential velocity such that
$$(z_{1}^{0}(\al) - \al, z_2^{0},(\al), \Phi^0_\al(\alpha)) \in X_{r_0},$$
for some $r_0>0$, and  satisfying:
\begin{enumerate}
\item $\displaystyle \left(u^{0}\cdot\frac{(z_\al^0)^\bot}{|z_\al^0|}\right)(\al_1) < 0, \left(u^{0}\cdot\frac{(z_\al^0)^\bot}{|z_\al^0|}\right)(\alpha_2)<0,$
\item $\displaystyle \int_{\partial \Omega}u^{0}\cdot\frac{(z_\al^0)^\bot}{|z_\al^0|} ds = \int_{\T} u^{0}(\al)\cdot (z_\al^0)^\bot d\al = 0$.
\end{enumerate}
Then there exist a finite time $T > 0$, $0<r<r_0$, a time-varying curve $\tilde{z}(\al,t)$ satisfying:
\begin{enumerate}
\item $P^{-1}(\tilde{z}_{1}(\al,t)) - \al, P^{-1}(\tilde{z}_{2}(\al,t))$ are $2\pi$-periodic,
\item $P^{-1}(\tilde{z}(\al,t))$ satisfies the arc-chord condition for all $t \in (0,T]$,
\end{enumerate}
and $\tilde{u}(\al,t)$ with
 $$(\tilde{z}_{1}(\al,t) - \al, \tilde{z}_2(\al,t), \tilde{\Phi}_\alpha(\alpha,t)(\al,t))\in C([0,T],  X_{r})$$
 which provides a solution of the water waves equations \eqref{1paco} with $\tilde{z}^0(\al)=P(z^0(\al))$ and $\tilde{u}(\al,0)\cdot(\tilde{z}_\alpha)^\perp(\alpha,0)=u^0(\al)\cdot(\tilde{z}^0)^\perp_\alpha(\alpha)$.
\end{thm}

 The main tool in the proof is to use an abstract  Cauchy-Kowaleswki theorem from \cite{Nirenberg:abstract-cauchy-kowalewski} and \cite{Nishida:theorem-Nirenberg} as in for example \cite{Castro-Cordoba-Fefferman-Gancedo-LopezFernandez:rayleigh-taylor-breakdown}.

\subsection{Local existence for initial data in Sobolev spaces in the tilde domain}

We will take the following $\tilde{c}(\al,t)$:

 \begin{align*}\tilde{c}(\al,t) & = \frac{\al+\pi}{2\pi}\int_{-\pi}^{\pi}(Q^{2}BR(\tilde{z},\tilde{\omega}))_\beta(\beta,t)\cdot\frac{\tilde{z}_{\beta}(\beta,t)}{|\tilde{z}_{\beta}(\beta,t)|^{2}}d\beta \\
& - \int_{-\pi}^{\al}(Q^{2}BR(\tilde{z},\tilde{\omega}))_\beta(\beta,t)\cdot\frac{\tilde{z}_{\beta}(\beta,t)}{|\tilde{z}_{\beta}(\beta,t)|^{2}}d\beta
\end{align*}

We will also define an auxiliary function $\tilde{\varphi}(\al,t)$ analogous to the one introduced in \cite{Beale-Hou-Lowengrub:growth-rates-linearized} (for the linear case) and \cite{Ambrose-Masmoudi:zero-surface-tension-waterwaves} (nonlinear case) which helps us to bound several of the terms that appear:
\begin{equation}\label{fvar}
\tilde{\varphi}(\al,t) = \frac{Q^2(\al,t)\tilde{\omega}(\al,t)}{2|\tilde{z}_{\al}(\al,t)|} - \tilde{c}(\al,t)|\tilde{z}_{\al}(\al,t)|.
\end{equation}

\begin{thm}
\label{localexistencetilde}
Let $z^{0}(\al)$ be a splash curve such that $z_{1}^{0}(\al) - \al, z_2^{0}(\al) \in H^{4}(\T)$. Let $u^{0}(\al) \cdot (z^{0}_{\al})^{\perp}(\al) \in H^4(\T)$ satisfying:
\begin{enumerate}
\item $\displaystyle \left(u^{0}\cdot\frac{(z_\al^0)^\bot}{|z_\al^0|}\right)(\al_1) < 0, \left(u^{0}\cdot\frac{(z_\al^0)^\bot}{|z_\al^0|}\right)(\alpha_2)<0,$
\item $\displaystyle \int_{\partial \Omega}u^{0}\cdot\frac{(z_\al^0)^\bot}{|z_\al^0|} ds = \int_{\T} u^{0}(\al)\cdot (z_\al^0)^\bot d\al = 0$.
\end{enumerate}
Then there exist a finite time $T > 0$, a time-varying curve $\tilde{z}(\al,t)  \in C([0,T];H^{4})$ satisfying:
\begin{enumerate}
\item $P^{-1}(\tilde{z}_{1}(\al,t)) - \al, P^{-1}(\tilde{z}_{2}(\al,t))$ are $2\pi$-periodic,
\item $P^{-1}(\tilde{z}(\al,t))$ satisfies the arc-chord condition for all $t \in (0,T]$,
\end{enumerate}
and $\tilde{u}(\al,t)\in C([0,T];H^3(\T))$ which provides a solution of the water waves equations \eqref{tildastream} with $\tilde{z}^0(\al)=P(z^0(\al))$ and $\tilde{u}(\al,0)\cdot(\tilde{z}_\alpha)^\perp(\alpha,0)=u^0(\al)\cdot(\tilde{z}^0)^\perp_\alpha(\alpha)$.
\end{thm}

\begin{sket}

In the proof, for the sake of simplicity, we will drop the tildes from the notation.

The proof will use the properties of $c(\al,t)$ and $\varphi(\al,t)$ to get an extra cancellation to help us derive energy estimates. Moreover, this choice of $c$ will ensure that the length of the tangent vector of $z(\al,t)$ depends only on time.

Here we define the energy $E(t)$ by
\begin{align*}
E(t)=&\|z\|^2_{H^3}(t)+\int_{\T}\frac{Q^2\sigma_z}{|z_\al|^2}|\da^4 z|^2d\al+\|F(z)\|^2_{L^\infty}(t)\\
&+\|\omega\|^2_{H^{2}}(t)+
\|\varphi\|^2_{H^{3+\frac12}}(t)+\frac{|z_\al|^2}{m(Q^2\sigma_z)(t)}+\sum_{l=0}^4\frac{1}{m(q^l)(t)}
\end{align*}
where the $L^\infty$ norm of the function
$$
F(z)\equiv \frac{|\beta|}{|z(\al,t)-z(\al-\beta,t)|},\quad \al,\beta\in\T
$$
measures the arc-chord condition,
\begin{align}
\begin{split}\label{R-T}
 \sigma_{z}  \equiv& \left(BR_{t}(z,\omega) + \frac{\varphi}{|z_{\al}|}BR_{\al}(z,\omega)\right) \cdot z_{\al}^{\perp} + \frac{\omega}{2|z_{\al}|^{2}}\left(z_{\al t} + \frac{\varphi}{|z_{\al}|}z_{\al \al}\right) \cdot z_{\al}^{\perp} \\
& + Q\left|BR(z,\omega) + \frac{\omega}{2|z_{\al}|^{2}}z_{\al}\right|^{2}(\nabla Q)(z) \cdot z_{\al}^{\perp}
 + (\nabla P_{2}^{-1})(z) \cdot z_{\al}^{\perp}
\end{split}
\end{align}
is the Rayleigh-Taylor function,
$$
m(Q^2\sigma_z)(t)\equiv\min_{\al\in\T}Q^2(\al,t)\sigma_z(\al,t),
$$
and finally
$$
m(q^l)(t)\equiv\min_{\al\in\T}|z(\al,t)-q^l|
$$
for $l=0,...,4$. We proceed as in \cite{Cordoba-Cordoba-Gancedo:interface-water-waves-2d}: The bound for the operator $(I+J)^{-1}$, where $J\omega=2BR(z,\omega)\cdot z_\alpha$,
 and some rather routine estimates allow us to find
$$
\frac{d}{dt}\Big(\|z\|^2_{H^3}(t)+\|F(z)\|^2_{L^\infty}(t)+\|\omega\|^2_{H^{2}}(t)+
\|\varphi\|^2_{L^2}(t)+\frac{|z_\al|^2}{m(Q^2\sigma_z)(t)}+\sum_{l=0}^4\frac{1}{m(q^l)(t)}\Big)\leq C E^k(t),
$$
for $C$ and $k$ universal constants. Above we use that
$$
\|\da^4 z\|^2_{L^2}(t)=\int_{\T}\frac{Q^2\sigma_z|z_\al|^2}{Q^2\sigma_z|z_\al|^2}|\da^4z|^2d\al
\leq \frac{|z_\al|^2}{m(Q^2\sigma_z)}\int_{\T}\frac{Q^2\sigma_z}{|z_\al|^2}|\da^4z|^2d\al\leq E^2(t).
$$
Further we obtain
$$
\frac{d}{dt}\Big(\int_{\T}\frac{Q^2\sigma_z}{|z_\al|^2}|\da^4 z|^2d\al\Big)\leq CE^k(t)+S(t)
$$
where
$$
S(t)=2\int_{\T}\frac{Q^2\sigma_z}{|z_\al|^2}\da^4 z\cdot\frac{z^\bot_\al}{|z_\al|}H(\da^4\varphi)d\al.
$$
We use \eqref{fvar}, \eqref{eqomega} and \eqref{R-T} to get
$$
\da^3\varphi_{t} =-\frac{Q^2\sigma_z}{|z_\al|^2}\da^4z\cdot\frac{z_\al^\bot}{|z_\al|}+\mbox{``control"}
$$
where ``control" is given by lower order terms and unbounded terms (such us $\da^4\varphi$) that can be estimated with energy methods in terms of $E(t)$. Therefore, it allows us to get
$$
\frac{d}{dt}\|\Lambda^{1/2}\da^3\varphi\|^2_{L^2}(t)\leq CE^k(t)-S(t)
$$
which together with above inequalities yields
$$
\left|\frac{d}{dt}E(t)\right|\leq CE^k(t).
$$
Local existence follows using standard arguments with the apriori energy estimate.

\end{sket}

\section{Structural Stability}

Again, in this section, we will omit the tildes from the notation. This section is devoted to establish a stability result. It will allow us to conclude the following: if $(x,\gamma)$ approximately satisfies equation \eqref{CharlieFlat} , then near to $(x,\gamma)$ there exists an exact solution $(z,\omega)$. Below is the theorem.

\begin{thm}
\label{stabilitytheorem}
Let
$$ D(\al,t) \equiv z(\al,t) - x(\al,t), \quad d(\al,t) \equiv \omega(\al,t) - \gamma(\al,t), \quad \mathcal{D}(\al,t) \equiv \varphi(\al,t) - \zeta(\al,t)$$
where $(x,\gamma,\zeta)$ are the solutions of
\begin{equation}
\label{CharlieFlat}
\left\{
\begin{array}{cl}
x_t & = Q^2(x)BR(x,\gamma) +  b x_{\al} + f\\
b & = \underbrace{\frac{\al + \pi}{2\pi}\int_{-\pi}^{\pi}(Q^2
BR(x,\gamma))_{\al}\frac{x_\al}{|x_{\al}|^{2}} d\al - \int_{-\pi}^{\al}(Q^2 BR(x,\gamma))_{\beta}\frac{x_{\al}}{|x_{\al}|^{2}}d\beta}_{b_s} \\
& + \underbrace{\frac{\al + \pi}{2\pi}\int_{-\pi}^{\pi}f_{\al}\frac{x_\al}{|x_{\al}|^{2}} d\al - \int_{-\pi}^{\al}f_{\beta}\frac{x_{\beta}}{|x_{\beta}|^{2}}d\beta}_{b_e}\\
\gamma_t & + 2BR_{t}(x,\gamma) \cdot x_{\al} = - (Q^2(x))_{\al}|BR(x,\gamma)|^{2} + 2bBR_{\al}(x,\gamma) \cdot x_{\al} + (b\gamma)_{\al} - \left(\frac{Q^2(x)\gamma^2}{4|x_{\al}|^{2}}\right)_{\al} \\
& \qquad \qquad \qquad - 2(P^{-1}_{2}(x))_{\al} + g \\
\zeta(\al,t) & = \frac{Q^2_x(\al,t)\gamma(\al,t)}{2|x_{\al}(\al,t)|} - b_s(\al,t)|x_{\al}(\al,t)|,
\end{array}
\right.
\end{equation}
$(z,\omega,\varphi)$ are the solutions of \eqref{CharlieFlat} with $f \equiv g \equiv 0$ and $\mathcal{E}$ the following norm for the difference
$$\mathcal{E}(t) \equiv \left(\|D\|^{2}_{H^{3}} + \int_{-\pi}^{\pi}\frac{Q^2\sigma_{z}}{|z_{\al}|^{2}}|\partial^{4}_{\al}D|^{2}d\alpha + \|d\|^{2}_{H^{2}} + \|\mathcal{D}\|^{2}_{H^{3+\frac{1}{2}}}\right).$$
Then we have that
$$\left|\frac{d}{dt}\mathcal{E}(t)\right|\leq \mathcal{C}(t)(\mathcal{E}(t)+\delta(t))$$
where $$\mathcal{C}(t)= \mathcal{C}(E(t),\|x\|_{H^{5+\frac12}}(t),\|\gamma\|_{H^{4+\frac12}}(t),
\|\zeta\|_{H^{4+\frac12}}(t),\|F(x)\|_{L^\infty}(t))$$ and $$\delta(t)=(\|f\|_{H^{5+\frac12}}(t)+\|g\|_{H^{3+\frac12}}(t))^k + (\|f\|_{H^{5+\frac12}}(t)+\|g\|_{H^{3+\frac12}}(t))^2, $$ with $k$ big enough. Here $E(t)$ is defined in the proof of Theorem \ref{localexistencetilde}.
\end{thm}

\begin{sket}
The equation for $(x,\gamma,\zeta)$ is the same as the one for $(z,\omega,\varphi)$ but for $f$ and $g$. The function $b$ is chosen in such a way that $|x_\al|$ only depends on time. Then it allows us to get the following estimates:
$$
\frac{d}{dt}\big(\|D\|^{2}_{H^{3}}+\|d\|^{2}_{H^{2}}+\|\mathcal{D}\|^{2}_{L^2}\big)\leq \mathcal{C}(t)(\mathcal{E}(t)+\delta(t)).
$$
Further we obtain
$$
\frac{d}{dt}\Big(\int_{\T}\frac{Q^2\sigma_z}{|z_\al|^2}|\da^4 D|^2d\al\Big)\leq \mathcal{C}(t)(\mathcal{E}(t)+\delta(t))+\mathcal{S}(t)
$$
with
$$
\mathcal{S}(t)=2\int_{\T}\frac{Q^2\sigma_z}{|z_\al|^2}\da^4 D\cdot\frac{z^\bot_\al}{|z_\al|}H(\da^4\mathcal{D})d\al.
$$
For $\mathcal{D}$ one finds that
$$
\da^3\mathcal{D}_{t} =-\frac{Q^2\sigma_z}{|z_\al|^2}\da^4D\cdot\frac{z_\al^\bot}{|z_\al|}+\mbox{``control"}
$$
where ``control" denotes terms which can be estimated by $\mathcal{C}(t)(\mathcal{E}(t)+\delta(t))$, which yields
$$
\frac{d}{dt}\|\Lambda^{1/2}\da^3\mathcal{D}\|^2_{L^2}(t)\leq
\mathcal{C}(t)(\mathcal{E}(t)+\delta(t))-\mathcal{S}(t).
$$
Then the desired estimate follows.
\end{sket}

\section{Numerical results}

In order to illustrate the splash singularity, several numerical simulations were performed. The simulations were done following the scheme proposed by Beale, Hou and Lowengrub \cite{Beale-Hou-Lowengrub:convergence-boundary-integral} adapted to the equations on the tilde domain (i.e. taking into account the impact of $Q$ on the equation). Instead of having an evolution equation for $\omega$, they introduce a velocity potential $\phi$ and study its evolution through time subject to the constraint imposed by being a potential. This is the set of equations \eqref{equations-hou-tilde}. The initial data on the non-tilde domain was given by:

$$ z^{0}_{1}(\al) = \al + \frac{1}{4}\left(-\frac{3\pi}{2} - 1.9\right)\sin(\al)+\frac{1}{2}\sin(2\al)+\frac{1}{4}\left(\frac{\pi}{2} - 1.9\right)\sin(3\al)$$
$$ z^{0}_{2}(\al) = \frac{1}{10}\cos(\al) - \frac{3}{10}\cos(2\al) + \frac{1}{10}\cos(3\al)$$

Note that $z\left(\frac{\pi}{2}\right) = z\left(-\frac{\pi}{2}\right)$ (splash). Instead of prescribing an initial condition for $\omega$, we prescribed the normal component of the velocity to ensure a more controlled direction of the fluid. From that we got the initial $\omega(\al,0)$ using the following relations. Let $\psi$ be such that $\nabla^{\perp}\psi = v$ and $\Psi(\al)$ its restriction to the interface. Recall that we can transform the initial condition on the normal component of the velocity into an initial condition on the tangential component by applying the transformations described in section 2. The initial normal velocity is then prescribed by setting
$$ u^{0}_{n}(\al)|z_{\al}(\al)| = \Psi_{\al}(\al) = 3 \cdot \cos(\al) - 3.4 \cdot \cos(2\al) + \cos(3\al) + 0.2\cos(4\al).$$

\begin{figure}\centering
\includegraphics[scale=0.2]{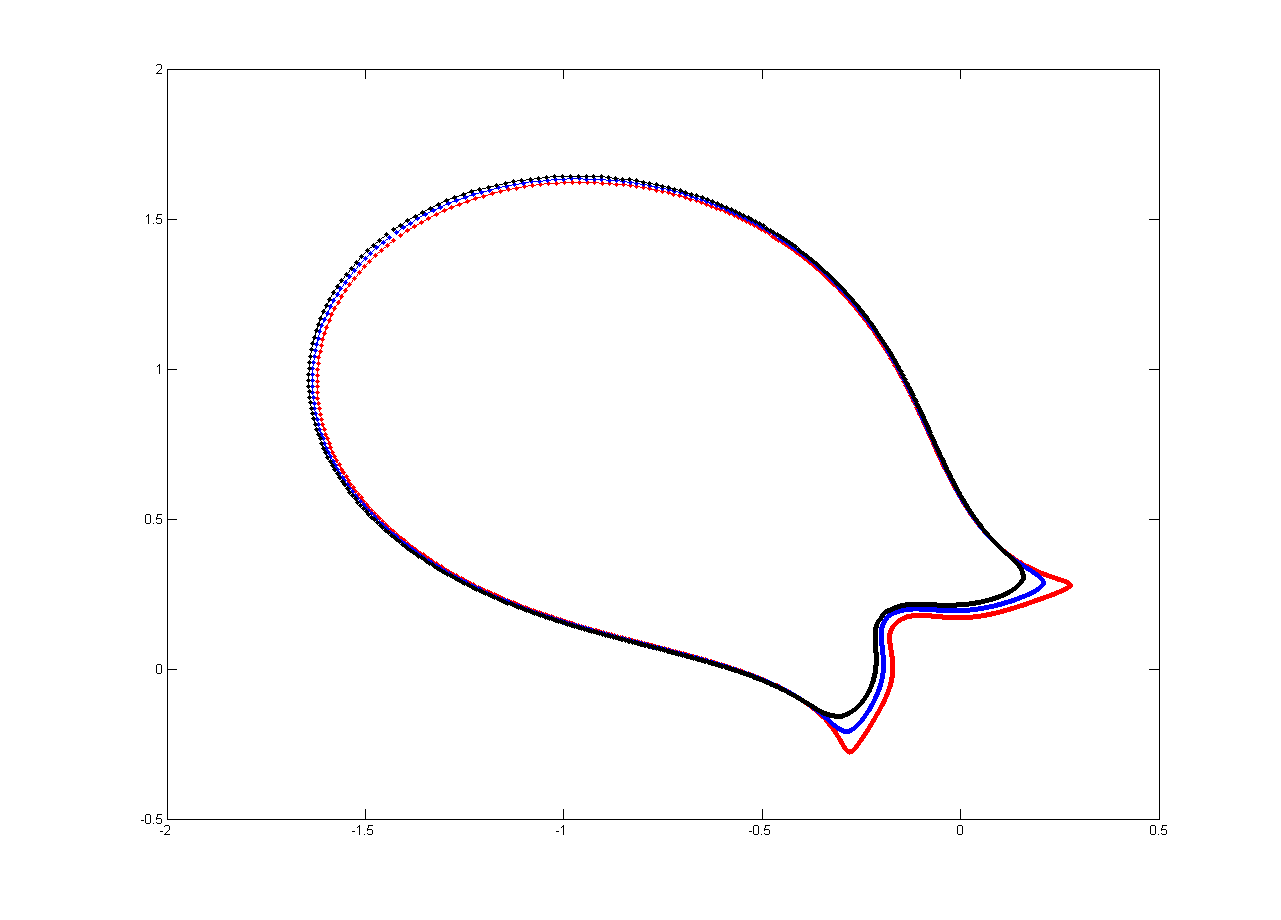}
\caption{Tilde domain at times $t = 0$ (Red - splash. The splash points are located on the open segment joining the points $(\sqrt{2}/2,\sqrt{2}/2)$ and $(-\sqrt{2}/2,-\sqrt{2}/2)$, approximately at $(\pm 0.27284156, \pm 0.27284156)$), $t = 4 \cdot 10^{-3}$ (Blue - turning) and $t = 7 \cdot 10^{-3}$ (Black - graph)}
\label{PictureSplash}
\includegraphics[scale=0.2]{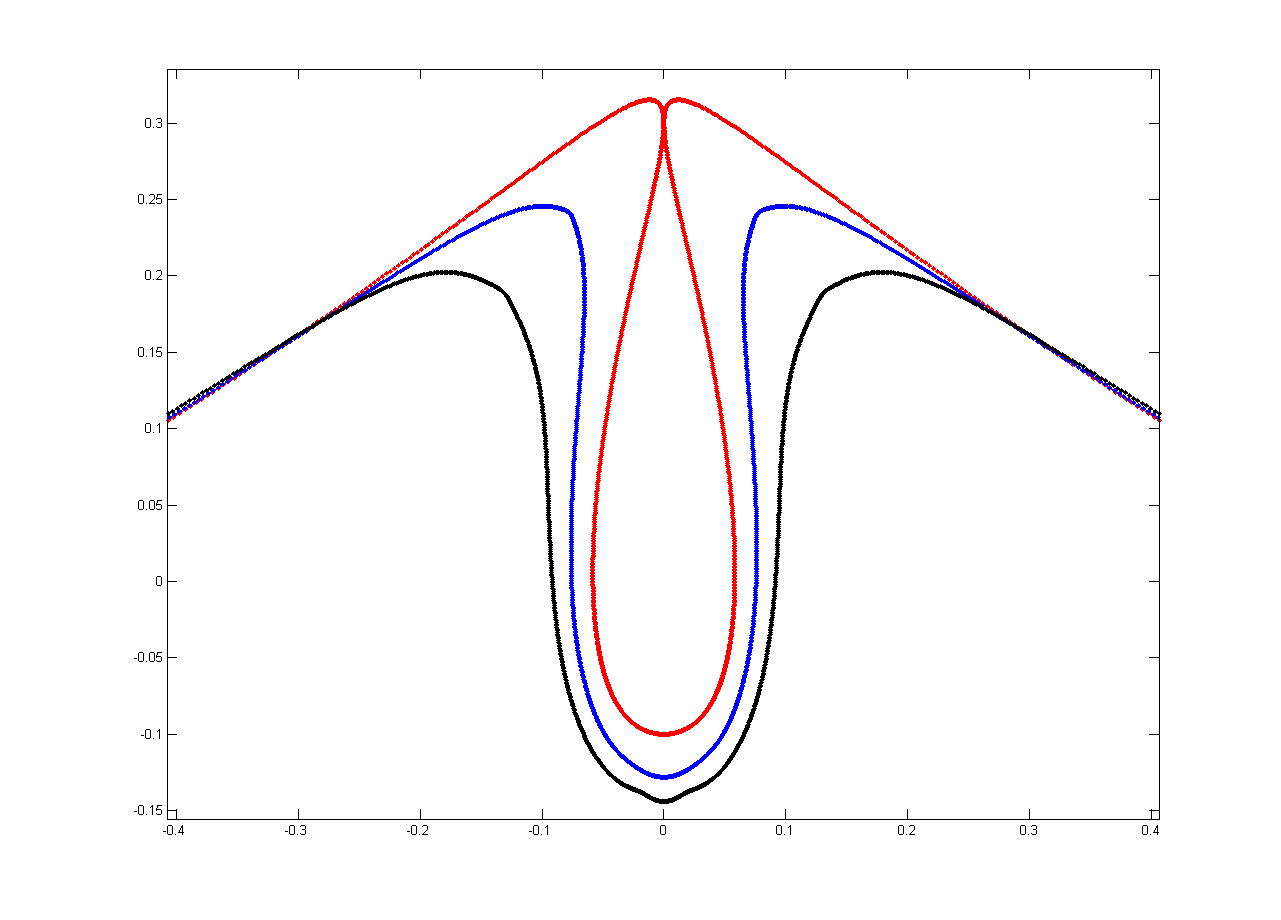}
\caption{Zoom of the splash singularity at times $t = 0$ (Red - splash), $t = 4 \cdot 10^{-3}$ (Blue - turning) and $t = 7 \cdot 10^{-3}$ (Black - graph)}
\label{PictureNonTildaZoom}
\end{figure}

The simulations were done using a spatial mesh of $N = 2048$ nodes and a time step $\Delta t = 10^{-7}$. The time direction was set to run backwards (from the splash to the graph) and the graph was obtained at approximately $T_{g} = 6.5 \cdot 10^{-3}$. Note that the normal component of the velocity $u^{0}_{n}\left(\pm \frac{\pi}{2}\right) > 0$ at the splash, which satisfies the hypotheses of Theorem \ref{localexistencenontilde} as we are running time backwards. Getting the potential of the initial condition from $\omega_0$ and $z_{0}$ and the transformation of all the initial data to the tilde domain is a trivial computation, taking care to choose the appropriate branch of the square root.  See figures \ref{PictureSplash} and \ref{PictureNonTildaZoom}.

\section{Further research}

We would like to exhibit a water-wave solution whose interface starts as an $H^4$-smooth graph at time zero, and ends in a splash at time $T$. We sketch a few ideas that may lead to a rigorous computer-assisted proof of the existence of such a solution. We will work in the tilde domain.

A simulation as in Section 5 leads to an approximate solution $x(\al,t)$, $\gamma(\al,t)$ with the desired properties. Thus $x(\cdot,t)$ describes a graph when $t = 0$ and a splash when $t = T$. Moreover, we believe that equations similar to \eqref{CharlieFlat} hold, with very small $f$ and $g$.

We may suppose that $x$ and $\gamma$ are known piecewise-polynomial functions on $[0,2\pi] \times [0,T]$. Using interval arithmetic \cite{Moore-Bierbaum:methods-applications-interval-analysis}, one can compute rigorous upper bounds for appropriate Sobolev norms of $f$ and $g$. We hope that these upper bounds will be very small.

Next, we solve the water-wave equations (\ref{em}-\ref{cEuler}) for $z,\omega$, starting at time $T$, and proceeding backwards in time. We take our initial $(z,\omega)$ at time $T$ to be a splash, very close to $(x(\cdot,T),\gamma(\cdot,T))$ in a high Sobolev norm. \footnote{Since $x$ and $\gamma$ are piecewise polynomials, we should not assume that $|\da x(\al,t)|$ is a function of $t$ alone. Hence, we cannot take $(z,\omega) = (x,\gamma)$ at time $T$.}

We want to compare the exact solution $(z,\omega)$ with the approximate solution $(x,\gamma)$, using the quantity $\mathcal{E}(t)$ as in Theorem \ref{stabilitytheorem}. Since $(z,\omega)$ and $(x,\gamma)$ are very close at time $T$, we will be able to show easily that

\begin{equation}
\label{Charlie!}
\mathcal{E}(T) < \varepsilon_1, \text{ for a very small, computable constant $\varepsilon_1$.}
\end{equation}

Moreover, the functions $x$ and $\gamma$ are known; and we also know upper bounds for Sobolev norms of $f$ and $g$. Therefore, the ideas in the proof of Theorem \ref{stabilitytheorem}, together with interval arithmetic, should lead to a rigorous proof of the differential inequality

\begin{equation}
\label{Charlie!!}
\left|\frac{d}{dt}\mathcal{E}(t)\right| \leq C_1 \mathcal{E}(t) + \varepsilon_2,
\end{equation}

where $C_1$ and $\varepsilon_2$ are computable constants.

We hope that $\varepsilon_2$ is very small, since $f$ and $g$ have small Sobolev norms; and we hope that $C_1$ won't be too big.

Once we establish \eqref{Charlie!} and \eqref{Charlie!!}, we will then know that our water wave solution $(z,\omega)$ exists for all time $t \in [0,T]$, and that $\mathcal{E}(0) < \varepsilon_3$ for a (hopefully small) computable constant $\varepsilon_3$.

From the definition of $\mathcal{E}(t)$, we will then easily deduce that $z(\cdot,0) - x(\cdot,0)$ has norm at most $\varepsilon_4$ in $H^{4}(\mathbb{R}/2\pi\mathbb{Z})$, for a computable constant $\varepsilon_4$.

If $\varepsilon_4$ is small enough, this in turn implies that the interface $z(\cdot,0)$ is an $H^{4}$-smooth graph. Thus $(z,\omega)$ is an exact solution of the water-wave equation, whose interface is an $H^4$-smooth graph at time 0, and a splash at time $T$.

We hope that a proof along these lines can be made to work.

\subsection*{{\bf Acknowledgements}}

\smallskip

 AC, DC, FG and JGS were partially supported by the grant {\sc MTM2008-03754} of the MCINN (Spain) and
the grant StG-203138CDSIF  of the ERC. CF was partially supported by
NSF grant DMS-0901040. FG was
partially supported by NSF grant DMS-0901810.

\bibliographystyle{abbrv}
\bibliography{references}

\begin{tabular}{ll}
\textbf{Angel Castro} &  \\
{\small Instituto de Ciencias Matem\'aticas} & \\
{\small Consejo Superior de Investigaciones Cient\'ificas} &\\
{\small C/ Nicol\'{a}s Cabrera, 13-15} & \\
{\small Campus Cantoblanco UAM, 28049 Madrid} & \\
{\small Email: angel\underline{  }castro@icmat.es} & \\
   & \\
\textbf{Diego C\'ordoba} &  \textbf{Charles Fefferman}\\
{\small Instituto de Ciencias Matem\'aticas} & {\small Department of Mathematics}\\
{\small Consejo Superior de Investigaciones Cient\'ificas} & {\small Princeton University}\\
{\small C/ Nicol\'{a}s Cabrera, 13-15} & {\small 1102 Fine Hall, Washington Rd, }\\
{\small Campus Cantoblanco UAM, 28049 Madrid} & {\small Princeton, NJ 08544, USA}\\
{\small Email: dcg@icmat.es} & {\small Email: cf@math.princeton.edu}\\
 & \\
\textbf{Francisco Gancedo} &  \textbf{Javier G\'omez-Serrano}\\
{\small Department of Mathematics} & {\small Instituto de Ciencias Matem\'aticas}\\
{\small University of Chicago} & {\small Consejo Superior de Investigaciones Cient\'ificas}\\
{\small 5734 University Avenue,} & {\small C/ Nicol\'{a}s Cabrera, 13-15} \\
{\small Chicago, IL 60637, USA}  & {\small Campus Cantoblanco UAM, 28049 Madrid} \\
{\small Email: fgancedo@math.uchicago.edu} & {\small Email: javier.gomez@icmat.es}\\
\end{tabular}

\end{document}